\documentclass[12pt,reqno]{amsart}
\setbox0=\hbox{$+$}
\newdimen\plusheight
\plusheight=\ht0
\def\+{\;\lower\plusheight\hbox{$+$}\;}

\setbox0=\hbox{$-$}
\newdimen\minusheight
\minusheight=\ht0
\def\-{\;\lower\minusheight\hbox{$-$}\;}

\setbox0=\hbox{$\cdots$}
\newdimen\cdotsheight
\cdotsheight=\plusheight
\def\cds{\lower\cdotsheight\hbox{$\cdots$}}

\renewcommand{\(}{\left\(}
\renewcommand{\)}{\right\)}

\def\NI{\noindent}

\renewcommand{\pmod}[1]{\,(\textup{mod}\,#1)}
\numberwithin{equation}{section}
 \theoremstyle{plain}
\newtheorem{theorem}{Theorem}[section]

\begin{document}
\footnotetext[2]{Corresponding author.}
\begin{center}{\bf \large General Congruences Modulo 5 and 7 for Colour Partitions}\end{center}\vskip
 5mm
\begin{center}\vskip 3mm
\centerline{\bf Nipen Saikia$^{1, 2}$ and Chayanika Boruah$^3$}\vskip3mm

{\it $^1$Department of Mathematics, Rajiv Gandhi University, \\Rono Hills, Doimukh, Arunachal Pradesh, India-791112.\\ E. Mail: nipennak@yahoo.com}\\

{\it $^3$Department of Mathematics, University of Science and Technology, \\ Ri-Bhoi, Meghalaya-793101, India.\\ E. Mail: cboruah123@gmail.com}

 \end{center}\vskip3mm

 \NI{\bf Abstract:}  For any positive integers $n$ and  $r$, let $p_r(n)$ denotes the number of partitions of $n$ where each part has $r$ distinct colours. Many authors studied the partition function $p_r(n)$ for  particular values of $r$. In this paper, we prove some general congruences modulo $5$ and $7$ for the colour partition function $p_r(n)$ by considering some general values of $r$. To prove the congruences  we employ some $q$-series identities which is also in the spirit of Ramanujan. \vskip
 3mm

 \NI {\bf Keywords and Phrases:} colour partition; $q$-series; congruence.\vskip
 3mm

 \NI{\bf Mathematics Subject Classifications:} 11P82; 11P83.\vskip
 3mm

 \section{Introduction}A partition of a positive integer $n$ is a non-increasing sequence of positive integers, called parts, whose sum equals $n$. For example, $n=3$ has three partitions, namely, 
  	$$3,\quad 2+1,\quad 1+1+1.$$ If $p(n)$ denote the number of partitions of $n$, then $p(3)=3$. The generating function for  $p(n)$ is given by
  	\begin{equation}\label{pngen}\sum_{n=0}^\infty p(n)q^n=\dfrac{1}{(q; q)_\infty},
  	\end{equation} where, here and throughout the paper 
  	$$(a;q)_{\infty}=\prod_{n=0}^{\infty}(1-aq^n).$$
  	Ramanujan \cite{sr} established following beautiful congruences for  $p(n)$:
  	$$p(5n+4) \equiv0\pmod5,$$
  	$$p(7n+5)\equiv0 \pmod7,$$ and
  	$$p(11n+6)\equiv0\pmod {11}.$$ In this paper we are concerned with colour partitions of  positive integer $n$.  	
  	 A part in a partition of $n$ has $r$ colours if there are
  	  $r$ copies of each part available and all of them are viewed as distinct objects. For example, if each part in the partition of 3 has two colours, say red and green, then the number of two colour partitions of $3$ is 10, namely  $$3_r,\quad 3_g, \quad2_r+1_r, \quad2_r+1_g,\quad 2_g+1_g, \quad 2_g+1_r,$$$$ 1_r+1_r+1_r,\quad 1_g+1_g+1_g,\quad 1_r+1_g+1_g,\quad 1_r+1_r+1_g.$$ Thus, number of 2-colour partitions of 3 is 10.
  	  The generating function of $r$-colour partitions of any positive integer $n$ is connected to the  general partition function $p_r(n)$ introduced by Ramanujan in a letter to Hardy \cite{bc1} and is given by  
  	  \begin{equation}\label{gf} \sum_{n=0}^{\infty}p_r(n)q^n=\dfrac{1}{(q;q)^r_{\infty}}.\end{equation}For $r=1$, $p_1(n)$ is the usual unrestricted partition function $p(n)$  defined in \eqref{pngen}.  	  
  	  If $r$ is negative, then \begin{equation}\label{mr}p_r(n)=\left(p_r(n, e)-p_r(n, o)\right),\end{equation} where $p_r(n, e)$ (resp. $p_r(n, o)$) is the number of partitions of $n$ with even (resp. odd) number of distinct parts  and each part have $r$ colours. For example, if $n=5$ and $r=-1$ then $p_{-1}(5, e)=2$ with relevant partitions $4+1$ and $3+2$, and $p_{-1}(5, o)=1$ with the relevant partition 5. Thus, $p_{-1}(5)=2-1=1$. Similarly, we see that $p_{-2}(3)=4-2=2$. The case $r=-1$ in \eqref{mr} is the famous Euler's pentagonal number theorem.  	  
  	  Ramanujan \cite{bc1} showed that, if $\lambda$ is a positive integer and  $\overline{w}$ is a prime of the form $6\lambda-1$, then  
  	  \begin{equation}
  	  p_{-4}\left(n\overline{w}-\frac{(\overline{w}+1)}{6}\right)\equiv0\pmod{\overline{w}}.
  	  \end{equation} 
  	  Ramanathan \cite{kg}, Atkin \cite{ao}, and Ono \cite{ko} investigated the partition function for some negative values of $r$. Recently, Saikia and Chetry \cite{anali} proved some infinite families of congruences modulo 7 for the partiton function $p_r(n)$ for negative values of $r$. 
  	  
  	   For positive values of $r$,  $p_r(n)$  counts the number of $r$-colour partitions of a positive integer $n$.  
 Gandhi \cite{JM} studied the colour partition function $p_r(n)$ for some particular values of $r$ and found some  Ramanujan-type congruences for  certain values of $r$.  For example, he proved that 
$$p_2(5n+3)\equiv 0\pmod 5 \quad \mbox{and} \quad
p_8(11n+4)\equiv 0\pmod{11}.$$ Newman \cite{ab} also found some congruences for colour partition. Baruah and Ojah \cite{NK} proved some congruences for $p_{3}(n)$ modulo some powers of 3. Recently, Hirschhorm \cite{Hirs} found congruences for $p_3(n)$ modulo higher powers of 3. 

In this paper, we prove some general congruences modulo 5 and 7 for the $r$-colour partition function $p_r(n)$  for some general values of $r$. To prove our congruences we will employ some $q-$series identities which is also  in the spirit of Ramanujan.  We list our congruences moduli 5 and 7 in Theorems \ref{mod5} and \ref{mod7}, respectively below: 

\begin{theorem}\label{mod5}For any non-negative integer $k$, we have
\begin{align}\hspace{-3.5cm}(i)&\quad
p_{5k+1}(5n+4)\equiv 0~\pmod5. \notag\\
\hspace{-3.5cm}(ii)&\quad p_{5k+2}(5n+i)\equiv 0~\pmod5, \quad for \quad  i= 2, 3,  4.\notag\\
  \hspace{-3.5cm}(iii)& \quad p_{5k+4}(5n+i)\equiv 0~\pmod5,\quad for \quad  i=3, 4.\notag\\
\hspace{-3.5cm}(iv)&\quad
 	p_{25k+3}(25n+22)\equiv 0~\pmod5.\notag\\
\hspace{-3.5cm} (v)&\quad 
 	 	p_{25k+4}(25n+21)\equiv 0~\pmod5.\notag
 \end{align}
\end{theorem}

\begin{theorem}\label{mod7} For any non-negative integer $k$, we have
\begin{align}\hspace{-2.5cm}(i)&\quad p_{7k+1}(7n+5)\equiv 0~\pmod7.\notag\\
\hspace{-2.5cm}(ii)&\quad p_{7k+4}(7n+j)\equiv 0~\pmod7,\quad for\quad j=2, 4, 5, 6.\notag\\
\hspace{-2.5cm}(iii)&\quad p_{7k+6}(7n+j)\equiv 0~\pmod7,\quad for\quad j=3, 4, 6.\notag\\
\hspace{-2.5cm}(iv)&\quad p_{49k+2}(49n+7j+3)\equiv 0~\pmod7, \quad for\quad j=2, 4, 5, 6.\notag\\
\hspace{-2.5cm}(v)&\quad p_{49k+3}(49n+7j+1)\equiv 0~\pmod7,\quad for\quad j=2, 4, 5, 6.\notag\\
\hspace{-2.5cm}(vi)& \quad p_{49k+5}(49n+39)\equiv 0~\pmod7.\notag
\end{align}
\end{theorem}

 \section{Preliminaries}  Ramanujan \cite{collected} stated that
 \begin{equation}\label{l1}(q;q)_{\infty}=(q^{25};q^{25})_{\infty}(F^{-1}(q^5)-q-q^2F(q^5)),
 \end{equation}
 where $F(q):=q^{-1/5}R(q)$ and $R(q)$ is the Rogers-Ramanujan continued fraction given by $$R(q):=q^{1/5}\dfrac{(q^2;q^5)_{\infty}(q^3;q^5)_{\infty}}{(q;q^5)_{\infty}(q^4;q^5)}=\frac{q^{1/5}}{1}_{+}\frac{q}{1}_{+}\frac{q^2}{1}_{+}\frac{q^3}{1}_{+
 	\cdots}, \hspace{.4cm}\vert q\vert<1.$$
 From \eqref{l1}, it is easy to see that
 \begin{equation}\label{rq2}(q; q)_\infty^2=(q^{25};q^{25})^2_{\infty}\Big(F^{-2}(q^5)-2qF^{-1}(q^5)-q^2+2q^3F(q^5)+q^4F^2(q^5)\Big),
 \end{equation}
 \begin{equation}\label{rq3}
 (q;q)_\infty^3=(q^{25};q^{25})^3_{\infty}\Big(F^{-3}(q^5)-3qF^{-2}(q^5)+5q^3-3q^5F^2(q^5)-q^6F^3(q^5)\Big),
 \end{equation}and
 $$(q; q)_\infty^4 =(q^{25};q^{25})^4_{\infty}\Big(F^{-4}(q^5)-4qF^{-3}(q^5)+2q^2F^{-2}(q^5)+8q^3F^{-1}(q^5)-5q^4$$\begin{equation}\label{rq4}\hspace{2.5cm}-8q^5F(q^5)+2q^6F^2(q^5)+4q^7F^3(q^5)+q^8F^4(q^5)\Big).
 \end{equation} 

Again, by \cite[p. 303, Entry 17(v)]{bcb3}, we have
\begin{equation}\label{l2}
\left(q;q\right)_\infty=(q^{49};q^{49})_\infty\left(A(q^7)-qB(q^7)-q^2+q^5C(q^7)\right),
\end{equation}
where
$$A(q^7)=\frac{f(-q^{14},-q^{35})}{f(-q^7,-q^{42})},~ B(q^7)=\frac{f(-q^{21},-q^{28})}{f(-q^{14},-q^{35})},~ C(q^7)=\frac{f(-q^7,-q^{42})}{f(-q^{21},-q^{28})},$$ and $$f(a,b)=\sum_{n=-\infty}^\infty a^{n(n+1)/2}b^{n(n-1)/2},\quad |ab|<1.$$
Squaring \eqref{l2}, we find that
$$\hspace{-.8cm}(q;q)^2_\infty=(q^{49};q^{49})^2_\infty\Big(\left(A(q^7)^2-2q^7C(q^7)\right)-2qA(q^7)B(q^7)+q^2\left(B(q^7)^2-2A(q^7)\right)$$\begin{equation}\label{eqj2}\hspace{2.8cm}+q^3\left(2B(q^7)+q^7C(q^7)^2\right)+q^4+2q^5A(q^7)C(q^7)-2q^6B(q^7)C(q^7)\Big).
\end{equation} 
Also, from \cite[p. 39, Entry 24(ii)]{bcb3} we note that
\begin{equation}\label{eqj1}
(q;q)_\infty^3=\sum_{n=0}^{\infty}(-1)^n(2n+1)q^{n(n+1)/2}.
\end{equation}
From \eqref{eqj1}, it follows that
$$\hspace{-1.1cm}(q;q)^3_\infty=J_{0}(q^7)-qJ_{1}(q^7)+q^3J_{3}(q^7)-7q^6J_{6}(q^7)$$
\begin{equation}\label{eqj3}\equiv{J_{0}(q^7)}-qJ_{1}(q^7)+q^3J_{3}(q^7)\pmod7
\end{equation} and 
$$\hspace{1.3cm}
(q;q)^6_\infty\equiv{J_{0}(q^7)^2}-2qJ_{0}(q^7)J_{1}(q^7)+q^2J_{1}(q^7)^2+2q^3J_{0}(q^7)J_{3}(q^7)$$
\begin{equation}\label{eqj4}\hspace{.8cm}-2q^4J_{1}(q^7)J_{3}(q^7)+q^6J_{3}(q^7)^2\pmod7,
\end{equation}
where $J_{0}, J_{1}, J_{3}$, and $J_{6}$ are series with integral powers of $q^7$. 

 From \cite[Lemma 3.12]{garvan} we note that, if 	
\begin{equation}\label{xi} \xi=\dfrac{(q;q)_{\infty}}{q^2(q^{49};q^{49})_{\infty}},\quad T_7=\dfrac{(q^7;q^7)^4_{\infty}}{q^7(q^{49};q^{49})^4_{\infty}}\end{equation} and $H_7$ is an operator which acts on a series of powers of  $q$ and picks out those terms in which the power of $q$ is congruent to $0$ modulo $7$, then
 \begin{equation}\label{l3} 
 	 H_7(\xi^4)=-4T_7-7 \quad \mbox{and} \quad 
 H_7(\xi^5)=10T_7+49.
 	\end{equation} 
In addition to the above $q-$identities, we will also need the following congruence which follows from the binomial theorem (or see \cite[Lemma 2.4]{ns}): For any prime $p$, we have
\begin{equation}\label{qp}(q^p;q^p)_{\infty}\equiv(q;q)_{\infty}^p~\pmod p.\end{equation}

\section{Proof of Theorem \ref{mod5}}

\noindent{\it Proof of $(i)$}:
 Setting $r=5k+1$ in \eqref{gf}, we obtain
\begin{equation}\label{qs1}\sum_{n=0}^{\infty}p_{5k+1}(n)q^n=\frac{1}{(q;q)_{\infty}^{5k+1}},
\end{equation}
Using \eqref{qp} in \eqref{qs1}, we obtain
\begin{equation}\label{qs2}\sum_{n=0}^{\infty}p_{5k+1}(n)q^n\equiv\frac{(q;q)^4_{\infty}}{(q^5;q^5)_{\infty}^{k+1}}~\pmod5, \end{equation}
Employing \eqref{rq4} in \eqref{qs2} and then extracting terms involving $q^{5n+4}$, dividing by $q^4$, and replacing $q^5$ by $q$, we arrive at the desired result.

\vskip3mm
\noindent{\it Proof of $(ii)$}: Setting $r=5k+2$ in \eqref{gf}, we obtain	
\begin{equation}\label{qs4}\sum_{n=0}^{\infty}p_{5k+2}(n)q^n=\frac{1}{(q;q)_{\infty}^{5k+2}},
\end{equation}	
Using \eqref{qp} in \eqref{qs4}, we obtain
\begin{equation}\label{qs5}\sum_{n=0}^{\infty}p_{5k+2}(n)q^n\equiv\frac{(q;q)^3_{\infty}}{(q^5;q^5)_{\infty}^{k+1}}~\pmod5, \end{equation}
Employing \eqref{rq3} in \eqref{qs5} and extracting terms involving $q^{5n+i}$ for $i=2, 3, 4$, we arrive  at the desired result.

\vskip3mm
\noindent{\it Proof of $(iii)$}: Setting $r=5k+4$ in \eqref{gf}, we obtain \begin{equation}\label{qs7}\sum_{n=0}^{\infty}p_{5k+4}(n)q^n=\frac{1}{(q;q)_{\infty}^{5k+4}},
\end{equation}	Using \eqref{qp} in \eqref{qs7}, we obtain
\begin{equation}\label{qs8}\sum_{n=0}^{\infty}p_{5k+4}(n)q^n\equiv\frac{(q;q)_{\infty}}{(q^5;q^5)_{\infty}^{k+1}}~\pmod5, \end{equation}
Using \eqref{l1} in \eqref{qs8} and extracting terms containing $q^{5n+i}$ for $i=3, 4$, we complete the proof.

 \vskip3mm
 \noindent{\it Proof of $(iv)$}:
 Setting $r=25k+3$ in \eqref{gf}, we obtain	\begin{equation}\label{qs12}\sum_{n=0}^{\infty}p_{25k+3}(n)q^n=\frac{1}{(q;q)_{\infty}^{25k+3}},\end{equation}
Using \eqref{qp} in \eqref{qs12}, we obtain
\begin{equation}\label{qs14}\sum_{n=0}^{\infty}p_{25k+3}(n)q^n\equiv\frac{(q;q)^2_{\infty}}{(q^{25};q^{25})^k_{\infty}(q^5;q^5)_{\infty}}~\pmod5,
\end{equation}
 Employing \eqref{rq2} in \eqref{qs14} and then extracting terms involving $q^{5n+2},$ dividing by $q^2$, and replacing $q^5$ by $q$, we obtain
\begin{equation}\label{qs15}\sum_{n=0}^{\infty}p_{25k+3}(5n+2)q^n\equiv\frac{4}{(q;q)_{\infty}(q^5;q^5)_{\infty}^{k-2}}~\pmod5,\end{equation}
Simplyfing \eqref{qs15} by using \eqref{qp}, we obtain
\begin{equation}\label{qs16}\sum_{n=0}^{\infty}p_{25k+3}(5n+2)q^n\equiv\frac{4(q;q)^4_{\infty}}{(q^5;q^5)_{\infty}^{k-1}}~\pmod5,\end{equation}
Employing \eqref{rq4} in \eqref{qs16} and extracting terms involving $q^{5n+4}$, we complete the proof.

\vskip3mm
\noindent{\it Proof of $(v)$}:
 Setting $r=25k+4$ in \eqref{gf}, we obtain	\begin{equation}\label{qs17}\sum_{n=0}^{\infty}p_{25k+4}(n)q^n=\frac{1}{(q;q)_{\infty}^{25k+4}},\end{equation}
Using  \eqref{qp} in \eqref{qs17}, we obtain
\begin{equation}\label{qs19}\sum_{n=0}^{\infty}p_{25k+4}(n)q^n\equiv\frac{(q;q)_{\infty}}{(q^{25};q^{25})^k_{\infty}(q^5;q^5)_{\infty}}~\pmod5,
\end{equation}
Employing  \eqref{l1} in \eqref{qs19} and then extracting terms involving $q^{5n+1},$ dividing by $q$, and replacing $q^5$ by $q,$  we have
\begin{equation}\label{qs20}\sum_{n=0}^{\infty}p_{25k+4}(5n+1)q^n\equiv\frac{4(q;q)^4_{\infty}}{(q^5;q^5)^k_{\infty}}~\pmod5,\end{equation}
Again, employing \eqref{rq4} in \eqref{qs20} and extracting terms involving $q^{5n+4},$ we arrive at the desired result.

\section{Proof of Theorem \ref{mod7}}
\noindent{\it Proof of $(i)$}:
 Setting $r=7k+1$ in \eqref{gf}, we obtain
\begin{equation}\label{qs21}\sum_{n=0}^{\infty}p_{7k+1}(n)q^n=\frac{1}{(q;q)_{\infty}^{7k+1}},
\end{equation}
Using \eqref{qp} in \eqref{qs21}, we obtain
\begin{equation}\label{qs22}\sum_{n=0}^{\infty}p_{7k+1}(n)q^n\equiv\frac{(q;q)^6_{\infty}}{(q^7;q^7)_{\infty}^{k+1}}~\pmod7, \end{equation}
Employing \eqref{eqj4} in \eqref{qs22} and extracting the terms involving $q^{7n+5}$, we arrive at the desired result.

\vskip3mm
\noindent{\it Proof of $(ii)$}:
Setting $r=7k+4$ in \eqref{gf}, we have
\begin{equation}\label{qs26}\sum_{n=0}^{\infty}p_{7k+4}(n)q^n=\frac{1}{(q;q)_{\infty}^{7k+4}},	\end{equation}
Using \eqref{qp} in \eqref{qs26}, we obtain
\begin{equation}\label{qs27}\sum_{n=0}^{\infty}p_{7k+4}(n)q^n\equiv\frac{(q;q)^3_{\infty}}{(q^7;q^7)_{\infty}^{k+1}}~\pmod7, \end{equation}
 Employing \eqref{eqj3} in \eqref{qs27} and extracting the terms involving $q^{7n+j}$ for $j=2, 4, 5, 6$,  we complete the proof.

\vskip3mm
\noindent{\it Proof of $(iii)$}:
Setting $r=7k+6$ in \eqref{gf}, we obtain
\begin{equation}\label{qs31}\sum_{n=0}^{\infty}p_{7k+6}(n)q^n=\frac{1}{(q;q)_{\infty}^{7k+6}},\end{equation}
Using \eqref{qp} in \eqref{qs31}, we obtain
\begin{equation}\label{qs32}\sum_{n=0}^{\infty}p_{7k+6}(n)q^n\equiv\frac{(q;q)_{\infty}}{(q^7;q^7)_{\infty}^{k+1}}~\pmod7, \end{equation}
Employing \eqref{l2} in \eqref{qs32} and extracting terms involving in $q^{7n+j}$ for $j=3, 4, 6$, we complete the proof.

\vskip3mm
\noindent{\it Proof of $(iv)$}:
Setting $r=49k+2$ in \eqref{gf}, we find that
\begin{equation}\label{qs40}\sum_{n=0}^{\infty}p_{49k+2}(n)q^n=\frac{1}{(q;q)_{\infty}^{49k+2}}.\end{equation}	
Using  \eqref{qp} in \eqref{qs40}, we obtain
\begin{equation}\label{qs41}\sum_{n=0}^{\infty}p_{49k+2}(n)q^n\equiv\frac{(q;q)^5_{\infty}}{(q^{49};q^{49})_{\infty}^k(q^7;q^7)_{\infty}}~\pmod7. \end{equation}
Employing \eqref{xi} in \eqref{qs41}, we obtain
\begin{equation}\label{qs42}\sum_{n=0}^{\infty}p_{49k+2}(n)q^n\equiv\frac{\xi^5q^{10}}{(q^{49};q^{49})^{k-5}_{\infty}(q^7;q^7)_{\infty}}~\pmod7.\end{equation}
Extracting the terms involving $q^{7n+3}$ and using operator $H_7$ in \eqref{qs42}, we obtain
\begin{equation}\label{qs43}\sum_{n=0}^{\infty}p_{49k+2}(7n+3)q^{7n+3}\equiv\frac{H_7(\xi^5)q^{10}}{(q^{49};q^{49})^{k-5}_{\infty}(q^7;q^7)_{\infty}}~\pmod7,\end{equation}
Employing \eqref{l3} in  \eqref{qs43}, we obtain
\begin{equation}\label{qs44}\sum_{n=0}^{\infty}p_{49k+2}(7n+3)q^{7n+3}\equiv3q^3\frac{(q^7;q^7)^3_{\infty}}{(q^{49};q^{49})^{k-1}_{\infty}}~\pmod7,\end{equation}
Dividing \eqref{qs44} by $q^3$ and replacing $q^7$ by $q,$ we obtain
\begin{equation}\label{qs45}\sum_{n=0}^{\infty}p_{49k+2}(7n+3)q^n\equiv3\frac{(q;q)^3_{\infty}}{(q^7;q^7)^{k-1}_{\infty}}~\pmod7,\end{equation}
Employing  \eqref{eqj3} in \eqref{qs45} and extracting terms involving $q^{7n+j}$ for $j=2, 4, 5, 6$, we complete the proof.

\vskip3mm
\noindent{\it Proof of $(v)$}:
Setting $r=49k+3$ in \eqref{gf}, we obtain
\begin{equation}\label{qs50}\sum_{n=0}^{\infty}p_{49k+3}(n)q^n=\frac{1}{(q;q)_{\infty}^{49k+3}},\end{equation}
Using \eqref{qp} in \eqref{qs50}, we obtain
\begin{equation}\label{qs51}\sum_{n=0}^{\infty}p_{49k+3}(n)q^n\equiv\frac{(q;q)^4_{\infty}}{(q^{49};q^{49})_{\infty}^k(q^7;q^7)_{\infty}}~\pmod7, \end{equation}
Employing \eqref{xi} in \eqref{qs51}, we obtain
\begin{equation}\label{qs52}\sum_{n=0}^{\infty}p_{49k+3}(n)q^n\equiv\frac{\xi^4q^8}{(q^{49};q^{49})^{k-4}_{\infty}(q^7;q^7)_{\infty}}~\pmod7,\end{equation}
Extracting the terms involving $q^{7n+1}$ and using operator $H_7$ in \eqref{qs52}, we obtain
\begin{equation}\label{qs53}\sum_{n=0}^{\infty}p_{49k+3}(7n+1)q^{7n+1}\equiv\frac{H_7(\xi^4)q^8}{(q^{49};q^{49})^{k-4}_{\infty}(q^7;q^7)_{\infty}}~\pmod7,\end{equation}
Employing \eqref{l3} in \eqref{qs53}, we obtain
\begin{equation}\label{qs54}\sum_{n=0}^{\infty}p_{49k+3}(7n+1)q^{7n+1}\equiv 3q\frac{(q^7;q^7)^3_{\infty}}{(q^{49};q^{49})^k_{\infty}}~\pmod7,\end{equation}
Dividing \eqref{qs54} by $q$ and replacing $q^7$ by $q,$ we obtain
\begin{equation}\label{qs55}\sum_{n=0}^{\infty}p_{49k+3}(7n+1)q^n\equiv3\frac{(q;q)^3_{\infty}}{(q^7;q^7)^k_{\infty}}~\pmod7,\end{equation}
Employing \eqref{eqj3} in \eqref{qs55} and extracting terms involving $q^{7n+j}$ for $j=2, 4, 5, 6,$ we arive at the desired result.

\vskip3mm
\noindent{\it Proof of $(vi)$}:
 Setting $r=49k+5$ in \eqref{gf}, we find that
\begin{equation}\label{qs60}\sum_{n=0}^{\infty}p_{49k+5}(n)q^n=\frac{1}{(q;q)_{\infty}^{49k+5}},\end{equation}
Using \eqref{qp} in \eqref{qs60}, we obtain
\begin{equation}\label{qs61}\sum_{n=0}^{\infty}p_{49k+5}(n)q^n\equiv\frac{(q;q)^2_{\infty}}{(q^{49};q^{49})_{\infty}^k(q^7;q^7)_{\infty}}~\pmod7, \end{equation}
Employing  \eqref{eqj2} in \eqref{qs61}, extracting terms involving in $q^{7n+4}$, dividing by $q^4$ and replacing $q^7$ by $q$, we obtain
\begin{equation}\label{qs62}\sum_{n=0}^{\infty}p_{49k+5}(7n+4)q^n\equiv\frac{1}{(q^7;q^7)^{k-2}_{\infty}(q;q)_{\infty}}~\pmod7,\end{equation}
Using \eqref{qp} in \eqref{qs62}, we obtain
\begin{equation}\label{qs62e}\sum_{n=0}^{\infty}p_{49k+5}(7n+4)q^n\equiv\frac{(q;q)_{\infty}^6}{(q^7;q^7)^{k-1}_{\infty}}~\pmod7,\end{equation}
Employing \eqref{eqj4} in \eqref{qs62e} and extracting terms involving $q^{7n+5}$, we arrive at the desired result.

\section*{Compliance with Ethical Standards}

\textbf{Conflict of interest:} The author declares that there is no conflict of interest regarding the publication of this article.

\textbf{Human and animal rights:} The author declares that there is no research involving human participants and/
or animals in the contained of this paper.

 \end{document}